\documentclass[12pt,leqno]{amsart}
\usepackage{graphicx,epsfig,color}
\usepackage{amsmath,amsfonts,amssymb}
\usepackage{array}
\newcolumntype{K}{>{\centering\arraybackslash}m{1.5cm}}

\def\N{\hbox{\font\dubl=msbm10 scaled 1200 {\dubl N}}}

\def\d{{\rm{d}}}

\newtheorem{Theorem}{Theorem}

\title[Sums of Consecutive Prime Squares]{Sums of Consecutive Prime Squares}

\author[Janyarak Tongsomporn, Saeree Wananiyakul, J\"orn Steuding]{Janyarak Tongsomporn, Saeree Wananiyakul, J\"orn Steuding}

\date{January 2021}

\begin{document}

\begin{abstract}
We prove explicit bounds for the number of sums of consecutive prime squares below a given magnitude
\end{abstract}

\maketitle

\noindent {\footnotesize{
\smallskip

}}

\section{Motivation and the Main Result}

Early last year the authors learned that 2020 can be represented as a sum of squares of consecutive prime numbers, namely
$$
2020=17^2+19^2+23^2+29^2.
$$
It is a natural question to ask what the next year with this property will be. We shall show that such a representation is a rare event. 

Indeed, if ${\rm{scp}}(x)$ counts the number of sums of squares of consecutive primes below $x$, i.e.,
$$
{\rm{scp}}(x)=\sharp\left\{p_{n}^2+p_{n+1}^2+\ldots +p_{n-1+m}^2\leq x\,:\, m\in\N\right\},
$$
where $p_j$ denotes the $j$-th prime number in ascending order, then $\lim_{x\to\infty}{\rm{scp}}(x)/x=0$. The following theorem provides more precise bounds. 

\begin{Theorem}
We have
$$
2\,{x^{1/2}\over \log x}<\pi(\sqrt{x})\leq  {\rm{scp}}(x)<10.9558\,{x^{2/3}\over (\log x)^{4/3}}\,,
$$
where the inequality on the far left is valid for $x\geq 289$ and all those to the right for $x>1$.
\end{Theorem}

Here, as usual, $\pi(N)$ is counting the number of primes $p\leq N$ and explicit bounds for this prime counting function are the main tool for proving the inequalities above; we have chosen a recent paper \cite{dusart} by Pierre Dusart. The dear reader is invited to improve upon the bounds of our theorem; maybe it is even possible to prove an asymptotic formula for the number of sums of consecutive prime squares below a given magnitude. Note that we do not consider here the question whether or not an integer can have two or even more such representations or how many of these exist. 

Using a computer algebra package one can verify that the next sum of squares of consecutive primes is given by the expected suspect, namely
$$
2189=13^2+17^2+19^2+23^2+29^2.
$$
A list with all integers below $5000$ that can be written as a sum of consecutive prime squares can be found in the third and final section.

This year's prime factorization is $2021=43\cdot 47$ which is a product of two consecutive primes. Following our approach one can also discuss sums of products of two consecutive primes; the corresponding bounds should be close to those found here for sums of consecutive prime squares. 

\section{Proof of the Theorem}

It is convenient to define, for fixed $m\in\N$, the counting function for sums of $m$ consecutive prime squares, i.e.,
$$
{\rm{scp}}_m(x)=\sharp\left\{p_{n}^2+p_{n+1}^2+\ldots +p_{n-1+m}^2\leq x\right\}.
$$

For the lower bound we first observe that the number of squares of prime numbers $p^2$ below or equal to $x$ is given by $\pi(\sqrt{x})$. 

In the sequel we shall use the explicit bounds
\begin{equation}\label{dusa}
{N\over \log N}<\pi(N)< 1.2551\,{N\over \log N},
\end{equation}
where the left inequality is valid for $N\geq 17$ and the one on the right for $N>1$ (see Corollary 5.2 in \cite{dusart}); of course, the celebrated prime number theorem provides an asymptotic formula for $\pi(N)$ with main term $N/\log N$, however, for excluding the related error term for our analysis, we prefer the version above with the factor $1.2551$. The corresponding range for these inequalities (resp. the range for $x$ in our theorem) is also useful with respect to computer experiments. 

It follows from (\ref{dusa}) that
$$
{\rm{scp}}(x) \geq {\rm{scp}}_1(x)=\pi(\sqrt{x})>{\sqrt{x}\over {1\over 2}\log x},
$$
which is valid for $x\geq 17^2=289$. This proves the lower bound.

The reasoning for the upper bound is a little more advanced. First we note that for $n={\rm{scp}}(x)$ we have
$$
mp_n^2\leq p_{n}^2+p_{n+1}^2+\ldots +p_{n-1+m}^2\leq x< p_n^2+p_{n+1}^2+\ldots+p_{n-1+m}^2+p_{n+m}^2.
$$
Hence, by the inequality in (\ref{dusa}),  
\begin{equation}\label{up1}
{\rm{scp}}_m(x)=n\leq \pi(\sqrt{x/m})<1.2551\,{\sqrt{x/m}\over {1\over 2}\log(x/m)}\leq 2.5102\,{(x/m)^{1/2}\over \log x},
\end{equation}
which is valid for $x>m$, which, obviously, is no severe restriction (since the largest integer $\leq x$ is a trivial upper bound for the length of a sum of consecutive primes squares $\leq x$). 

To continue we shall next bound the length of possible sums of consecutive prime squares below $x$. For this purpose we shall use an old result due to Barkley Rosser \cite{rosser} which has been improved several times, in particular by Dusart \cite{dusart}, however, we prefer the more simple inequality
$$
p_n>n\log n,
$$
valid for all $n\in\N$; this lower bound is trivial for $n=1$. We shall use this so-called Rosser theorem for the sum of the squares of the {\it first} primes:
$$
p_{1}^2+p_{2}^2+\ldots +p_{M}^2>\sum_{2\leq n\leq M} (n\log n)^2.
$$
If we can show that the right hand side is larger than $x$, then the least sum of $M$ consecutive prime squares already exceeds the given magnitude. Assuming that this $M$ is the least positive integer with this property, this leads to a bound for $M$ depending on $x$. This estimate in combination with the previous one allows us to derive the upper bound of the theorem. Alternatively, one could also use partial summation here together with the prime number theorem, however, it is our intention to circumvent error terms. 

Obviously, for $M\geq 4$,
\begin{eqnarray*}
\sum_{2\leq n\leq M} (n\log n)^2&\geq& \sum_{\sqrt{M}\leq n\leq M} n^2(\log n)^2\\
&\geq& ({\textstyle{1\over 2}}\log M)^2\sum_{\sqrt{M}\leq n\leq M}n^2\geq {\textstyle{1\over 12}}M^3(\log M)^2,
\end{eqnarray*}
where we have used in the final step the well-known formula 
$$
1+2^2+3^2+\ldots+M^2={\textstyle{1\over 6}}\, M(M+1)(2M+1) 
$$
and some pen and paper. It thus follows that every sum of consecutive prime squares below $x$ has less than roughly $x^{1/3}$ summands. For a more precise bound we observe that substituting 
\begin{equation}\label{Mx}
M=\lfloor 108^{1/3}\, x^{1/3}(\log x)^{-2/3}\rfloor
\end{equation}
into the lower bound above yields a quantity slightly larger than $x$; here $\lfloor z\rfloor$ denotes the largest integer $\leq z$. 

To use this for an upper bound we first observe that (\ref{up1}) implies
$$
{\rm{scp}}(x)=\sum_{1\leq m\leq M}{\rm{scp}}_m(x)< 2.5102\, {x^{1/2}\over \log x}\sum_{1\leq m\leq M} m^{-1/2}.
$$
In general, we have, for $\alpha\in(0,1)$,
$$
\sum_{1\leq m\leq M} m^{-\alpha}< 1+\sum_{2\leq m\leq M}\int_{m-1}^m u^{-\alpha}\d u=1+\int_1^Mu^{-\alpha}\d u={M^{1-\alpha}-\alpha\over 1-\alpha}.
$$
This in combination with (\ref{Mx}) leads to
$$
{\rm{scp}}(x)< 5.0204\,{(xM)^{1/2}\over \log x}\leq {5.0204\cdot (108)^{1/6}}\,{x^{2/3}\over (\log x)^{4/3}}.
$$
This proves the upper bound of the theorem.

\section{Explicit Sums of Consecutive Prime Squares}

We conclude with a list of all integers below $5000$ that can be written as a sum of consecutive prime squares:
\begin{center}
\begin{tabular}{KKKKKKK}
4 & 9 & 13 & 25 & 34 & 38 & 49 \\
74 & 83 & 87 & 121 & 169 & 170 & 195 \\
204 & 208 & 289 & 290 & 339 & 361 & 364 \\
373 & 377 & 458 & 529 & 579 & 628 & 650 \\
653 & 662 & 666 & 819 & 841 & 890 & 940 \\
961 & 989 & 1014 & 1023 & 1027 & 1179 & 1348 \\
1369 & 1370 & 1469 & 1518 & 1543 & 1552 & 1556 \\
1681 & 1731 & 1802 & 1849 & 2020 & 2189 & 2209 \\
2310 & 2330 & 2331 & 2359 & 2384 & 2393 & 2397 \\
2692 & 2809 & 2981 & 3050 & 3150 & 3171 & 3271 \\
3320 & 3345 & 3354 & 3358 & 3481 & 3530 & 3700 \\
3721 & 4011 & 4058 & 4061 & 4350 & 4489 & 4519 \\
4640 & 4689 & 4714 & 4723 & 4727 & 4852 & 4899 
\end{tabular}
\end{center}


\vspace*{1cm}

\noindent Janyarak Tongsomporn, Saeree Wananiyakul, ${\mathcal W}$alailak University, School of Science, Nakhon Si Thammarat 80\,160, Thailand\hfill tjanyarak@gmail.com
\bigskip

\noindent J\"orn Steuding, Department of Mathematics, ${\mathcal W}$\"urzburg University, Am Hubland, 97\,218 W\"urzburg, Germany\hfill steuding@mathematik.uni-wuerzburg.de

\end{document}